\documentclass{article}
\usepackage[final]{neurips_2021}
\pdfoutput=1
\usepackage[utf8]{inputenc} 
\usepackage[T1]{fontenc}    
\usepackage{hyperref}       
\usepackage{booktabs}       
\usepackage{amsfonts}       
\usepackage{nicefrac}       
\usepackage{microtype}      
\usepackage{xcolor}         
\usepackage{amsmath}
\usepackage{epstopdf}
\usepackage{subfig}
\usepackage{comment}
\usepackage[pdftex]{graphicx}
\usepackage{bm}
\usepackage{wrapfig}
\setcitestyle{square,numbers}
\title{On automatic differentiation for the Mat\'ern covariance}

\author{%
Oana Marin\\
 Mathematics and Computer Science\\
  Argonne National Laboratory\\
   \texttt{oanam@anl.gov} \\
\And 
Christopher Geoga\\
Department of Statistics\\
Rutgers University\\
\texttt{christopher.geoga@rutgers.edu}\\
\And 
Michel Schanen\\
 Mathematics and Computer Science\\
  Argonne National Laboratory\\
   \texttt{mschanen@anl.gov} 
}

\begin{document}

\maketitle

\begin{abstract}
To target challenges in differentiable optimization we analyze and propose strategies for derivatives of the Mat\'ern kernel with respect to the smoothness parameter. This problem is of high interest in Gaussian processes modelling due to the lack of robust derivatives of the modified Bessel function of second kind with respect to order. 
In the current work we focus on newly identified series expansions for the modified Bessel function of second kind valid for complex orders. Using these expansions we obtain highly accurate results using the complex step method. Furthermore, we show that the evaluations using the recommended expansions are also more efficient than finite differences.
\end{abstract}

\section{Introduction}
The Mat\'ern kernel, widely used in Gaussian Processes (GP) statistical modeling, contains the modified Bessel function of the second kind, which is a special function. Special functions lack a formal definition and subsequently  robustness of implementation, even without considering automatic differentiation. We suggest the simple mnemonic that a function can be considered \emph{special} if it is not algebraic or lacks a predefined universal implementation. 
\begin{wrapfigure}[15]{r}{0.5\textwidth}
   \includegraphics[trim=60 240 90 300,clip,width=0.5\textwidth]{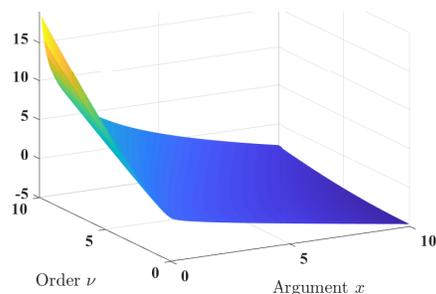}
   \centering
   \captionsetup{margin=1cm}
  \caption{The function $\log_{10}K_{\nu}(\bm{x})$ for a range of orders and arguments.\label{fig:besselk}}
\end{wrapfigure}  
The most problematic are functions such as the Bessel or the Anger function, or hyper-geometric functions, which typically require dedicated library implementations, e.g. \cite{amos1983a,cody1976}. These libraries, however, focus entirely on special function evaluations, may offer few implementation details, and may provide inaccurate results when automatically differentiated. Previous work on the automatic differentiation (AD) of certain special functions required the development of entirely new tools, e.g. \cite{charpentier2015, charpentier2018} to compute reliable derivatives. Very niche special functions, or less common derivatives, may not be available in automatic differentiation codes for special functions, which leaves a user with the single option of employing finite differences on special function evaluations.

Bessel functions are widely used in electromagnetics where they contribute to the kernel of the Helmholtz equation, or statistics where the modified Bessel functions enter the Mat\'ern covariance kernel.
The modified Bessel functions are the solution pair $I_{\nu}(x),\ K_{\nu}(x)$ of order $\nu$ of the following differential equation
$$
x^2 y'' + x y' - \left(x^2 + \nu^2 \right)y = 0\ ,
$$
known as the modified Bessel functions of first kind, $I_{\nu}(x)$, and second kind $K_{\nu}(x)$, illustrated in Fig.~\ref{fig:besselk}. According to the nature of the order $\nu$, which can be either real, complex, positive or not, a range of series expansion expressions are available for the solution pair $I_{\nu}(x),\ K_{\nu}(x)$. The range of the argument $x$ gives rise to different asymptotic expansions, and small argument expressions \cite{abramowitz65} may differ significantly from large argument expressions \cite{nemes2017}. 
Considering that most expressions are implemented as a truncated series expansion, these functions exhibit a high sensitivity to round-off errors, or platform and compiler variations, albeit only at computer precision level. 

Derivatives with respect to argument $x$ were successfully developed within AD  using the Fa\'a di Bruno generalization of Taylor series, see \cite{charpentier2015, charpentier2018}. However, derivatives with respect to order are considerably more intricate both mathematically and numerically. For the Bessel function of second kind we note in Fig.~\ref{fig:besselk} that the function has a pole at zero, and decays exponentially with increasing argument $x$ for any fixed order $\nu$. This implies that if we seek derivatives with respect to order $\nu$ at a sufficiently decayed $x$ we operate entirely at computer precision level.
Analytical expressions are available for certain derivatives of special functions, in particular for integer orders, e.g. $\nabla_{x} K_1(x)=-K_0(x)-K_1(x)/x$, but not all. Given that implementations depend on truncation tolerances, interlacing analytical expressions with numerically sensitive expressions can lead to numerical artifacts increasing with the differentiation order \footnote{\url{www.advanpix.com/2016/05/12/accuracy-of-bessel-functions-in-matlab/}}.
We recall that instruction level automatic differentiation (AD) of a power series truncated at a fixed level is not guaranteed to yield a derivative sufficiently converged at the same truncation level. To this end finite difference derivatives may provide superior accuracy in certain cases \cite{more2014}, although overall the error of finite difference calculations is highly prone to round-off errors due to cancellations. 

These considerations indicate that differentiable programming of special functions is first and foremost a mathematically unresolved issue. Differentiating through iterative processes, which also have a tolerance imposed convergence, has been considered within the AD literature, \cite{bruce1994,gilbert1992,griewank1993}, however, it has not been applied for series expansions. Furthermore, this approach may induce a computational overhead via extra computations if the differentiated series has very slow convergence. A notable strategy that alleviates the round-off error effects of finite differences, while preserving the computational complexity of forward mode automatic differentiation is the complex step method \cite{fike2012automatic}. Provided the intricacy of the mathematical expression in complex space the complex step method has been largely disregarded as a suitable approach for special functions. With the advent of neural network models, and their inherent automatic differentiation companion, special function neural network models (SFNN) were considered as a stand-in for series expansions \cite{marin2021}. Although SFNNs are a very promising approach, the lack of uniqueness of neural network models, and subsequently possible portability issues across platforms still prevents their use in production codes. 

In the current work, and motivated by statistical modelling problems, we focus on robust derivatives $\partial_{\nu}K_{\nu}(x)$ for the entire spectrum of real positive $\nu$ and $x$. A background on the motivation is expanded in Sec.~\ref{sec:apps}, while in Sec.~\ref{sec:numerics} we present and compare different differentiation approaches. Numerical results, outlining the trade-off between accuracy and efficiency of implementation are illustrated in Sec.~\ref{sec:results}. 

\section{Applications of special functions}\label{sec:apps}
Considering that special functions are mostly encountered as kernels of integral equations, which is traditionally an academic field, the computational cost of multi-precision computations at small scales has been acceptable.
The surge of statistical modelling, particularly for large data sets, revives the demand for robust implementations of special function derivatives.  

\paragraph{Statistical modelling: Mat\'ern Class Covariance Functions} 

Gaussian process (GP) models are ubiquitous throughout the physical and
numerical sciences. They are incredibly convenient in many ways, for example
providing linear conditional expectations and being specified entirely by their
first two moments. A primary problem in the study of GPs is to properly specify
the \emph{covariance function} $C_{\bm{\theta}}$, indexed by parameters
$\bm{\theta}$, that determines the covariance of a Gaussian random field $Z$ at
indices $\bm{x}$ and $\bm{x}'$ (like spatial locations), so that
\begin{equation*} 
  \text{Cov}(Z(\bm{x}), Z(\bm{x}')) := C_{\bm{\theta}}(\bm{x}, \bm{x}')\ ,
\end{equation*}
under that model.  Among other purposes, the covariance function is crucial for determining the
behavior of interpolants and forecasts, both with respect to predicted values
and inferred second-order information.

Arguably the most popular covariance function arising in 
variety of applications in the physical sciences is the \emph{Mat\'ern class} of
covariance functions \citep{matern1960}, given by
\begin{equation} \label{eq:matern}
  C_{\sigma, \rho, \nu}(\Vert \bm{x} - \bm{x}'\Vert) 
  = 
  \frac{\sigma^2 2^{1-\nu}}{\Gamma(\nu)}
  \left(\frac{ \sqrt{2 \nu} \Vert \bm{x} - \bm{x}'\Vert}{\rho} \right)^{\nu}
  \mathcal{K}_{\nu} \left(
  \frac{ \sqrt{2 \nu} \Vert \bm{x} - \bm{x}'\Vert}{\rho} 
  \right).
\end{equation}
The Mat\'ern class holds an advantage over
other positive-definite functions, such as the squared exponential
or rational quadratic, by providing complete control of the
degree of smoothness of $K$ at the origin. In turn this determines the number
of \emph{mean-square derivatives} sample paths will have. Further, in the
common sampling regime of \emph{fixed-domain asymptotics}, where measurements
are made more and more densely in a fixed spatial domain, the mean-square
differentiability of a process is one of the few quantities that is actually
resolved better as more data is added to a sample. For more details, we refer
readers to \cite{stein1999}, particularly the discussion of micro-ergodicity and
results on equivalence and orthogonality of Gaussian measures.

The optimization required to perform
maximum likelihood estimation for Gaussian processes is computationally  challenging,
and the required derivatives of $C$ with respect to $\nu$ are not easily
computed. The likelihood itself in the mean-zero case for data $\bm{y}$ is given
by
\begin{equation*} 
  -2 \ell(\bm{\theta}) := 
  \log \left\vert \bm{\Sigma}(\bm{\theta}) \right\vert 
  + 
  \bm{y}^T \bm{\Sigma}(\bm{\theta})^{-1} \bm{y},
\end{equation*}
where $\bm{\Sigma}(\bm{\theta})_{j,k} = C_{\bm{\theta}}(\bm{x}_j, \bm{x}_k)$,
and the derivative with respect to $\theta_j$ is given by
\begin{equation*} 
  -2 (\nabla \ell(\bm{\theta}))_j = 
  \text{tr}\left(
  \bm{\Sigma}(\bm{\theta})^{-1} 
  \frac{\partial}{\partial\theta_j} \bm{\Sigma}(\bm{\theta}) 
  \right)
  - 
  \bm{y}^T 
  \bm{\Sigma}(\bm{\theta})^{-1} 
  \frac{\partial}{\partial\theta_j}\bm{\Sigma}(\bm{\theta})  
  \bm{\Sigma}(\bm{\theta})^{-1}
  \bm{y}.
\end{equation*}
Thus, one requires the derivatives $\frac{\partial}{\partial\theta_j}
C_{\bm{\theta}}(\bm{x}_j, \bm{x}_k)$ to assemble the derivative matrices and
compute the gradient of the likelihood. Beyond the general difficulty of being
nonlinear, the likelihood surface for many models is nearly flat along level
surfaces of certain nonlinear functions of several parameters due to them having
similar interpolation properties (see \cite{stein2013} for an example plot of
such a likelihood surface). While it is clearly valuable to fit the
smoothness parameter to data, as opposed to the current practice of fixing it ahead of time, an effective optimization demands reliable derivatives 
 of $C_\nu(\bm{x})$.

\section{Numerical approaches}\label{sec:numerics}
Various implementations and recommendations can be found in the literature of special functions. The NIST compendium \cite{NIST} provides a comprehensive list of identities, including the derivative $\nabla_{\nu}K_{\nu}(x)$, see \cite[Eq.10.38.2]{NIST}, considered valid for any real $\nu$. 
As we will illustrate this formula is unreliable, and different tracks have to be explored. Encouraged by earlier work in the space of automatic differentiation, which considers numerical analysis aspects to automate robust implementations \cite{charpentier2015,charpentier2018,griewank1993,fike2012automatic} we identify two series expansions which can be easily evaluated using the complex step method. 
\paragraph{Complex step method:}
Owing to its simplicity and efficiency the complex step method is of great appeal to automatic differentiation, see \cite{fike2012automatic}. As a quick overview we show how the Taylor expansion in complex space yields a reliable derivative evaluation. Consider the expansion around $x_0$

$$
f(x_0+ih)=f(x_0)+ih\nabla_{x}f(x_0)-\frac{h^2}{2!}\nabla^2_{x}f(x_0)-i\frac{h^3}{3!}\nabla^3_{x}f(x_0)\ .
$$
Restricting the discussion to first order derivatives, although higher orders are also possible, we have
\begin{equation*}\label{eq:taylor}
    ih\nabla_{x}f(x_0)=f(x_0+ih)-f(x_0)+\mathcal O(h^2)\ ,
\end{equation*}
which by equating the purely imaginary parts, provides a derivative as
\begin{equation}\label{eq:1stord}
\nabla_{x}f(x_0)\approx\frac{1}{h}\Im{f(x_0+ih)} \ .
\end{equation}

Note that Eq.~\ref{eq:1stord} involves the evaluation of a complex number, however only one division by a small $h$, which eliminates the effect of cancellations crippling finite difference evaluations. The same Taylor expansions considerations are used in deriving finite difference schemes. The expression in Eq.~\ref{eq:1stord} will yield a method of accuracy $\mathcal O(h)$, which is identical to that of first order finite differences. However, by having removed the subtraction of two very similar terms, and obtained a single evaluation the accuracy order $\mathcal O(h)$ is assured to be free of round-off errors for $h$ down to computer precision. This is not the case for finite differences which become very sensitive to floating point errors at values below $h\approx 10^{-8}$. Moreover, the implementation of the complex step method can be accelerated at compiler level using dual-types since a single evaluation suffices for obtaining the derivative.

\paragraph{Small arguments:}
To obtain a robust expression in complex space we performed a set of 
 analytical manipulations and substitutions, using formulas in \cite{abramowitz65,nemes2017,NIST}. The main focus was to avoid cancellations or ratios of close to zero-valued functions. We identified a numerically robust formula for the modified Bessel function valid for $\nu$ complex given by
\begin{equation}\label{eq:KseriesCS}
K_{\nu}(x) = \sum_{k=0}^{\infty}\frac{1}{2k!} \bigg( \frac{x}{2}\bigg)^{2k} \bigg(\Gamma(\nu) \bigg( \frac{x}{2}\bigg)^{-\nu}  \frac{\Gamma(1-\nu)}{\Gamma(1+k-\nu)} +\Gamma(-\nu) \bigg(\frac{x}{2}\bigg)^\nu \frac{\Gamma(1+\nu)}{\Gamma(1+k+\nu)}	\bigg), \
\end{equation}
where $\Gamma(\cdot)$ denotes the complex Gamma function. 
This series expansion is, however, accurate only up to argument values $x\leq 10$. In the complex step method the derivative requires evaluations of complex order for a preset $h$, which can be arbitrarily low since cancellations are no longer a concern, and we simply have
$$
\nabla_{\nu}K_{\nu}(x)=\frac{1}{h}\Im{K_{\nu+ih}(x)}.
$$
To implement this infinite sum we accumulate terms for the summation index $k$ up to an imposed tolerance as a stopping criterion, implying the sum is truncated at a certain level $N$. 

 \textbf{Large arguments: } As $x>10$ we identify a different complex order compatible expression for $\mathcal{K}_\nu$ in terms of the confluent hypergeometric function $U$, also known as Tricomi's function
\begin{equation}\label{eq:kummer}
 \mathcal{K}_{\nu}\left(z\right)=\sqrt{\pi}e^{-z}\left(2z
      \right)^{\nu}U\left(\nu+\tfrac{1}{2},2\nu+1,2z\right),
\end{equation} 
which can be written in terms of the second linearly independent confluent hypergeometric function $M$, also now as Kummer's function:
\begin{equation*}
U(a,b,z)=\frac{\Gamma(1-b)}{\Gamma(a+1-b)}M(a,b,z)+\frac{\Gamma(b-1)}{\Gamma(a)}z^{1-b}M(a+1-b,2-b,z),
\end{equation*}
where $M$ can be computed via the generalized hyper-geometric series 
\begin{equation*}
 M(a,b,z)=\sum_{n=0}^\infty \frac {a^{(n)} z^n} {b^{(n)} n!}={}_1F_1(a;b;z),
\end{equation*}
 with $a^{(0)}=1$, and $a^{(n)}=a(a+1)(a+2)\cdots(a+n-1)$. Note that Eq.~\ref{eq:kummer} is also an infinite series, which requires truncation, and in practice $N\approx 10$ provided close to computer precision accuracy. 

\section{Numerical results}\label{sec:results}
\begin{figure}[!ht]
\centering \subfloat[Fixed small argument $x=3.94$.]
           {\includegraphics[trim=90 220 100 180,clip,width=0.47\textwidth]{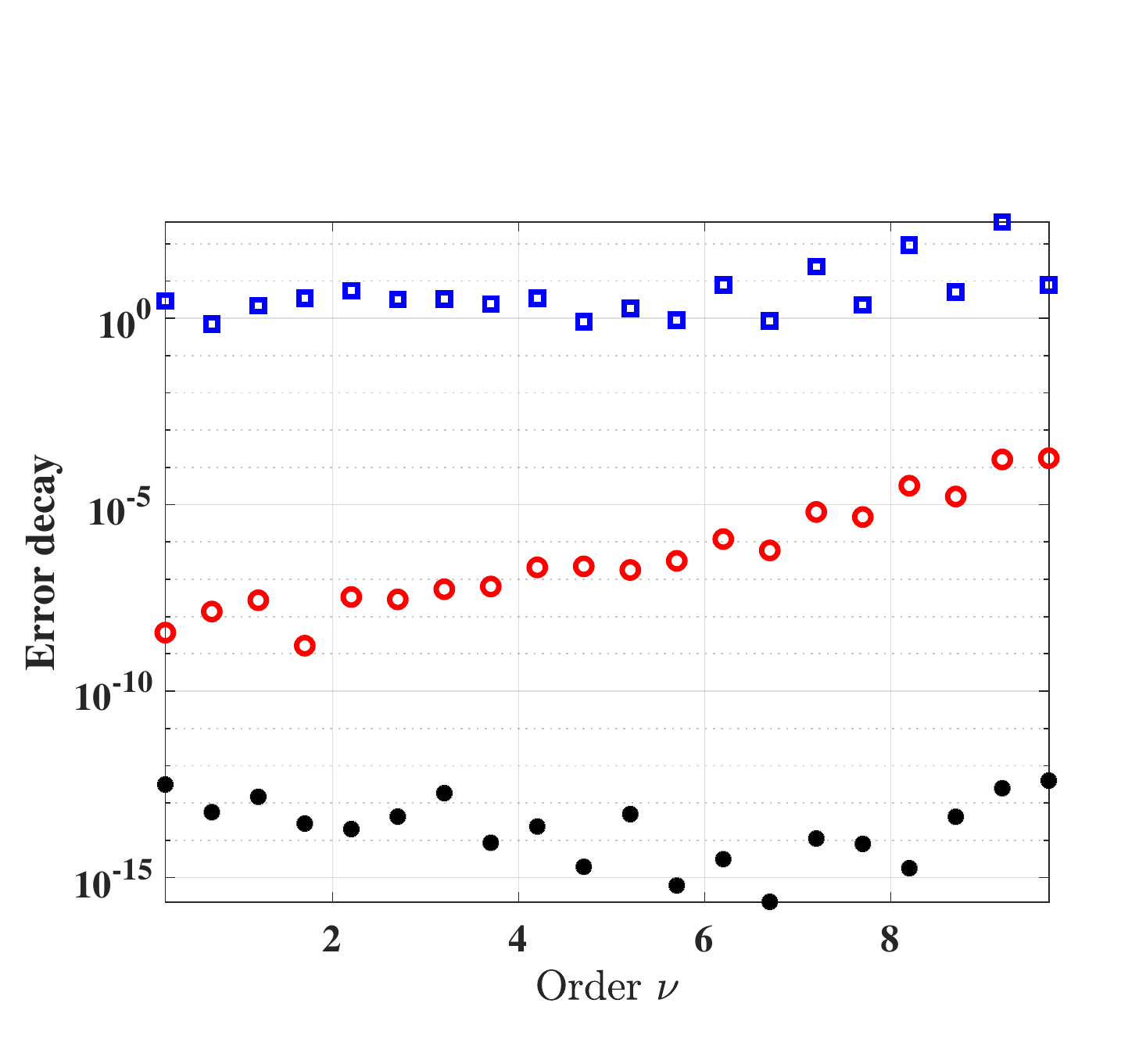}
\label{fig:err_small}}
\quad \subfloat[Fixed large argument $x=13.57$.]  {\includegraphics[trim=90 220 100 180,clip,width=0.47\textwidth]{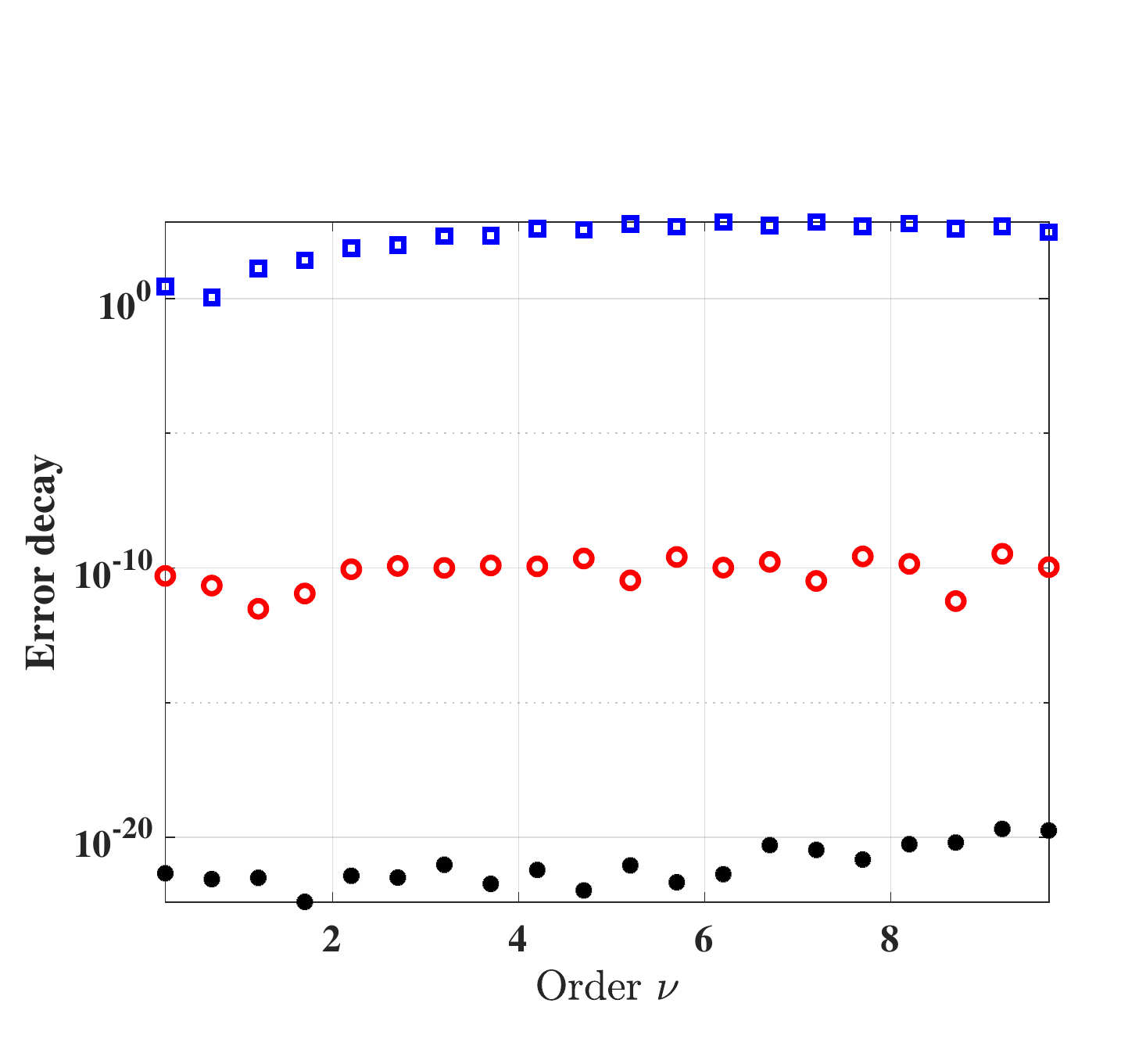}
\label{fig:err_large}}
\caption{Comparison of complex step (black markers), finite differences (red markers), naive implementation of \cite[Eq.:10.38.2]{NIST} (blue markers).}\label{fig:accuracy}
\end{figure}

In Fig.~\ref{fig:accuracy} we compare the accuracy of the complex step series implementation to a finite differences approach, and a direct naive  implementation of \cite[Eq.:10.38.2]{NIST}. The error is evaluated in absolute norm with  respect to a derivative computed in multi-precision to serve as ground truth. 
The series expansion derived in Eq.~\ref{eq:KseriesCS} is highly stable and accurate, and superior to finite differences for arguments $x<10$, as in Fig.~\ref{fig:err_small}. For increasingly large arguments $x>10$ both the finite difference implementation and the complex step method approach gain more correct digits, see Fig.~\ref{fig:err_large}, with the complex step being computer precision accurate.

Efficiency considerations are assessed in Tab.~\ref{tab:timings}, where we note the complex step method displays superior computational efficiency to second-order adaptive finite difference algorithm, using  \texttt{FiniteDifferences.jl}\footnote{L. White et al.. (2021). JuliaDiff/FiniteDifferences.jl: v0.12.18. \url{https://doi.org/10.5281/zenodo.5146583}}.  

\begin{table}[ht!]
    \centering
  \vspace{-0.5em}    
\begin{tabular}{|l|c|c|c|c|c|c|c|c|c|c|}
\hline
\multicolumn{1}{|c|}{$\bm{x}$} & \multicolumn{2}{c|}{0.05} & \multicolumn{2}{c|}{0.505}& \multicolumn{2}{c|}{1.05} & \multicolumn{2}{c|}{1.505}& \multicolumn{2}{c|}{2.05}\\
\cline{2-11}
\multicolumn{1}{|c|}{Method} & FD & CS & FD & CS & FD & CS & FD & CS & FD & CS \\
\hline
$\nu=0.25$&1.39 & 0.98&   1.52 &   1.32  &  1.63 &   1.60  & 1.74 &  1.70 &2.88  &   1.96 \\
$\nu=0.56$ &1.42 & 0.93 &  1.56 &  1.28 &  1.65 &  1.57 &  1.77 &   1.72  & 3.03 &  1.80 \\
$\nu=0.88$    &   1.37 &  0.93&  1.63 &  1.55 &  1.75 &   1.54 &   1.87 &  1.82  &  3.20 &   1.78 \\
$\nu=1.19$&    1.39 &  0.93&   1.54  &  1.45 &    1.77 &   1.76  &   1.90  & 1.82 &   3.00 &  1.79 \\
$\nu= 1.5$&    1.44 &   0.94&   1.58 &  1.46 &   1.67 &   1.66 &  1.79 &  1.85 &  2.11   &   2.05\\
\hline
\end{tabular} 
     \vspace{0.8em}
    \caption{Timing comparisons in microseconds for the computation of derivatives $\partial_{\nu} \mathcal{K}_{\nu}$; complex-step (CS) with fixed step $h=10^{-8}$ vs finite differences (FD).}
    \label{tab:timings}
      \vspace{-2em}
\end{table}

These expressions, although designed to be amenable to the complex step method, proved to have a superior performance also under \emph{black-box} AD, as we show in \cite{geogaAD2021}, where we extend this work. We have to specify the expressions provided here are not valid for edge-cases $\nu\in\mathbb Z$ or $\nu+\tfrac{1}{2}$, when we have to compute limit values. These limits are easy to derive and use in the complex step method, since they can be implemented as discrete values via \texttt{if} statements, however they cannot be easily ported to AD implementations which can assure higher efficiency. 

\section{Conclusions}\label{sec:conclusions}
The current work assesses state-of-the art strategies for differentiating Bessel functions with respect to order. The focus is on the modified Bessel function of second kind, $K_{\nu}(x)$, which is particularly relevant in statistical modelling as a component of the Mat\'ern covariance function. The function $K_{\nu}(x)$ is numerically cumbersome and prone to large errors under classical differentiation strategies. We identified two novel series expansions that are highly accurate and efficient to be used for derivative evaluations using the complex step method. 

\section*{Acknowledgments} The authors would like to thank Michael Stein (UChicago) for  inspiring and encouraging this work, as well as Lois Curfmann McInnes and Paul Hovland (Argonne National Laboratory) for engaging scientific discussions on automatic differentiation state-of-the-art strategies and neural network models. 
{\bf Government License.} The submitted manuscript has been created by
UChicago Argonne, LLC, Operator of Argonne National Laboratory
(``Argonne''). Argonne, a U.S. Department of Energy Office of Science
laboratory, is operated under Contract No. DE-AC02-06CH11357. The
U.S. Government retains for itself, and others acting on its behalf, a
paid-up nonexclusive, irrevocable worldwide license in said article to
reproduce, prepare derivative works, distribute copies to the public,
and perform publicly and display publicly, by or on behalf of the
Government.

\bibliography{spfun_short}

\end{document}